\documentclass{article}
\usepackage{amsmath,amsthm,amssymb}

\newcommand{\fc}{\mathfrak{c}}

\newcommand{\hm}{\mathfrak{hm}}

\newcommand{\ga}{\alpha}
\newcommand{\gb}{\beta}
\renewcommand{\gg}{\gamma}
\newcommand{\gd}{\delta}
\newcommand{\gw}{\omega}

\newcommand{\gS}{\Sigma}
\newcommand{\gs}{\sigma}
\newcommand{\oox}{\gw_{\gw_2+1}}
\newcommand{\fx}{\mathfrak{x}}

\newcommand{\R}{\mathbb{R}}

\newcommand{\cB}{\mathcal{B}}

\newcommand{\cov}{\mathtt{cov}}
\newcommand{\add}{\mathtt{add}}
\newcommand{\cof}{\mathtt{cof}}
\newcommand{\non}{\mathtt{non}}

\newtheorem{theorem}{Theorem}[section]
\newtheorem{lemma}[theorem]{Lemma}

\newtheorem{corollary}[theorem]{Corollary}
\newtheorem{fact}[theorem]{Fact}

\theoremstyle{definition}
\newtheorem{definition}[theorem]{Definition}

\title{Duality and the pcf theory\footnote{2000 AMS subject classification 03E17, 03E04}}
\author{
Saharon Shelah\thanks{Supported by The Israel Science 
Foundation founded by the Israel Academy of Sciences and Humanities. Publication 791.}
\\Hebrew University\\Rutgers University
\and
Jind{\v r}ich Zapletal\thanks{Partially supported by grants GA \v CR 201-00-1466, NSF DMS-0071437,
a CLAS UF research award and a visiting appointment at Hebrew University in February and March 2002.}
\\University of Florida}

\begin{document}

\maketitle
\begin{abstract}
We consider natural cardinal invariants $\hm_n$ and prove several duality theorems, saying roughly: if $I$
is a suitably definable ideal and provably $\cov(I)\geq\hm_n$, then $\non(I)$ is provably small. The proofs
integrate the determinacy theory, forcing and pcf theory.
\end{abstract}

\section{Introduction}

The authors of \cite{geschke:convexity} considered the following cardinal invariant of the continuum. Let $c:[2^\gw]^2\to 2$
be the partition defined by $c(x,y)=\Delta(x, y)\mod 2$ where $\Delta(x,y)=\min\{ m\in\gw:x(m)\neq y(m)\}$, and let
$J$ be the $\gs$-ideal $\gs$-generated by the homogeneous sets. It is not difficult to see that $J$ is a proper
$\gs$-ideal; in fact $c$ is the simplest continuous partition of $[2^\gw]^2$ with this property. \cite{geschke:convexity}
defined the cardinal invariant $\hm$ as the covering number of the ideal $J$. The invariant $\hm$ occurs naturaly in several 
contexts: it is the minimal number of Lipschitz-with-constant-1 functions from $2^\gw$ to itself such that
their graphs and graphs of their inverses cover the whole square $2^\gw\times 2^\gw$. It is also the minimal
number of convex subsets necessary to cover a certain closed subset of the Euclidean plane. See \cite{geschke:convexity}.

It is consistent with ZFC that $\hm<\fc$; not surprisingly \cite{z:detandci}, this is exactly what happens in the Sacks
model \cite{geschke:convexity}. On the other hand, $\hm$ is a very large cardinal invariant in that provably 
$\hm^+\geq\fc$. Here we prove a duality theorem similar to the one from \cite{z:detandci}.

\begin{theorem}
\label{t1}
Suppose that $I$ is an analytic $\gs$-ideal such that ZFC proves that
$\cov(I)\geq\hm$. Then ZFC proves that $\non(I)\leq\aleph_3$.
\end{theorem}

Here a $\gs$-ideal is analytic if there is a $\gS^1_1$ set in the plane such that the ideal is $\gs$-generated
by its vertical sections. Similarly, it is possible to define the class of projective ideals. The theorem
remains true for projective ideals if the theory ZFC is replaced by ZFC+``there are $\gw$ many Woodin cardinals''.

The theorem remains true if the pair $\cov, \non$ is replaced by any other dual pair of invariants, such as
$\non, \cov$ or $\add, \cof$. However, in order for the various complexity computations to come
out right, in these cases we must restrict ourselves to the class of Borel ideals.

The pattern persists to some variations of $\hm$. For every natural number $n>0$ let $c_n:[2^\gw]^2\to n$
be the function defined by $c_n(x, y)=\Delta(x,y)\mod n$ and let $J_n$ be the $\gs$-ideal on $2^\gw$
$\gs$-generated by the sets $X$ such that $c_n''X\neq n$. Define $J_1$ to be the ideal of countable sets.
Define $\hm_n$ to be the covering number of the ideal $J_n$. It is not difficult to show that
$\fc=\hm_1\geq\hm=\hm_2\geq\hm_3\geq\hm_4\geq\dots$ and $(\hm_n)^{n-1}\geq\fc$ for every number $n$. Again,
$\hm_n$ can be rewritten as the smallest number of functions in a certain class necessary to cover the cube $(2^\gw)^n$
with their graphs. We have

\begin{theorem}
\label{t2}
Let $n>0$. Suppose that $I$ is an analytic $\gs$-ideal such that ZFC proves that
$\cov(I)\geq\hm_n$. Then ZFC proves that $\non(I)\leq\aleph_{n+1}$.
\end{theorem}

\noindent In case of $n=1$ this improves the original duality theorem of \cite{z:detandci}, and the obtained bound
$\aleph_2$ is optimal. Finally,

\begin{theorem}
\label{t3}
Suppose that $I$ is an analytic ideal such that ZFC proves that
$\cov(I)\geq\min_n\hm_n$. Then ZFC proves that $\non(I)\leq\aleph_{\gw_2+1}$.
\end{theorem}

\noindent The remarks after Theorem~\ref{t1} remain in force for~\ref{t2} and~\ref{t3}.

The proofs of the above theorems integrate the effective descriptive theory, determinacy theory, forcing, and pcf theory.
This paper contains only the arguments that
are specific to the ideals $J_n$ in question. The general theory of definable proper forcing is encapsulated
into several facts, stated without proof. For the detailed development of this theory, the reader is referred
to the monograph \cite{z:detandci}.

The notation in the paper follows the set theoretic standard of \cite{jech:set}. Whenever $X$ is a set, $x\in X$
its element and $f:(2^\gw)^{X\setminus\{x\}}\to 2^\gw$ is a function, by the graph of the function $f$ we mean the
collection $\{\vec r\in (2^\gw)^X:\vec r(x)=f(\vec r\restriction X\setminus\{ x\})\}$.  For a $\gs$-ideal $I$ on the reals,
the cardinal $\cov(I)$ and $\non(I)$ denote the smallest size of a family of sets in the ideal covering
the whole real line, and the smallest size of an $I$-positive set respectively. If $B$ is an $I$-positive set,
then $I\restriction B$ is the ideal generated by $I$ and the complement of $B$.

\section{The forcings associated with $\hm_n$}

The following is a basic simple observation.

\begin{lemma}
\label{l1}
Fix a number $n>0$.
\begin{enumerate}
\item The ideal $J_n$ is $\gs$-generated by closed sets.
\item The ideal $J_n$ is nontrivial, i.e. $2^\gw\notin J_n$.
\end{enumerate}
\end{lemma}

As proved in \cite{geschke:coloring}, the partition $c_2$ is the minimal open partition of $[2^\gw]^2$ into two pieces such that the $\gs$-ideal
generated by the homogeneous sets is proper, in the sense that if $d$ is any other such a partition then there
is a Borel injection $i:2^\gw\to 2^\gw$ such that for all $x\neq y\in 2^\gw$ $d(i(x), i(y))=c_2(x,y)$. Similar
minimality result holds true for the partitions $c_n$.

\begin{proof}
The partitions $c_n$ are continuous. Therefore, if $X\subset 2^\gw$ is a generating set of the ideal $J_n$ 
(the image $c_n''[X]^2$ does not contain some number $m\in n$), then even its closure is such a generating set of the ideal
$J_n$ (its image under $c_n$ still does not contain that number $m$). The first item immediately
follows.

For the second item we will just show that the ideal $J_n$ is a subideal of the meager ideal, for every
number $n>0$. Fix the number $n$ and let $X\subset 2^\gw$ be a closed generating set of the ideal $J_n$ so that
the image $c_n''[X]^2$ leaves out some number $m\in n$. It will be enough to show that the set $X$ is nowhere
dense. And indeed, if $s\in 2^{<\gw}$ is a finite sequence representing some basic open set, prolong it
to obtain a sequence $t\supset s$ whose length is equal to $m$ modulo $n$. By the choice of the set $X$,
for one of the sequences $t^\smallfrown 0, t^\smallfrown 1$ the set $X$ contains no infinite binary
sequences extending it. Thus the set $X$ is nowhere dense as desired.
\end{proof}

Now look at the partial order $P_n$ of $J_n$-positive Borel sets ordered by inclusion. By Lemma 1.2 and 1.3 of 
\cite{z:closed}, this poset
adds a single real which falls out of all $J_n$-small sets and it is proper. In fact, this partial order has a natural
combinatorially simple dense subset. Call a tree $T\subset 2^{<\gw}$ $n$\emph{-fat} if it is nonempty and for every
node $s\in T$ and every number $m\in n$ there is a splitnode $t\supset s$ in the tree $T$ whose length
is equal to $m$ modulo $n$.

\begin{lemma}
\label{l5}
If $T$ is an $n$-fat tree then $[T]\notin J_n$. Moreover, if $A\subset 2^\gw$ is an analytic set, then either
it contains a subset of the form $[T]$ for some $n$-fat tree, or else it belongs to the ideal $J_n$.
\end{lemma}

\begin{proof}
The regularity property of analytic sets can be proved in several ways; we give a classical
determinacy argument following the proof of the perfect set theorem. 
Suppose that $A\subset 2^\gw$ is a set. Consider the game $G_n(A)$ where players
Adam and Eve alternate to play finite binary sequences $s_0, s_1, \dots$ and bits $b_0, b_1, \dots$ respectively such that
Adam's sequences $s_i$ form an extension increasing chain such that $s_i^\smallfrown b_i\subset s_{i+1}$
and the length of the sequence $s_i$ is equal to $i$ modulo $n$. Adam wins if $\bigcup_is_i\in A$.

First, Adam has a winning strategy in the game $G_n(A)$ if and only if the set $A$ contains all branches
of some $n$-fat tree. If $\gs$ is a winning strategy for Adam then the downward closure of the set of all
sequences that can arise in a play according to the strategy $\gs$ is an $n$-fat tree and all of its branches
belong to the set $A$. On the other hand, if $[T]\subset A$ is an $n$-fat tree then Adam can easily win the
game $G_n(A)$ by playing only splitnodes in the tree $T$.

Second, Eve has a winning strategy if and only if the set $A$ belongs to the ideal $J_n$. If $A\subset\bigcup_i X_i$
is in the ideal $J_n$, covered by countably many sets $X_i$ such that $i\mod n\notin c_n''[X_i]$, then Eve wins
by answering Adam's sequence $s_i$ in the $i$-th round with a bit $b_i$ such that no element of the
set $X_i$ begins with $s_i^\smallfrown b_i$. There is such a bit $b_i$ because no two elements of the set $X_i$
can first differ at the number $i\mod n$ which is the length of the sequence $s_i$. 
Of course, if Eve plays in this way she will win in the end,
since the real $\bigcup_is_i$ will fall out of all the sets $X_i$ and therefore out of the set $A$. On the other
hand, if Eve has a winning strategy $\gs$ then the set $A$ is in the ideal $J_n$. For every position $p$ consistent
with the strategy $\gs$ which ends after the round $i_p$ was completed with
some finite binary sequence $t_p$, let $S_p$ be the downward closure of the set
$\{t\in 2^{<\gw}:$ for no finite sequence $s$ and no bit $b$ it is the case that the play $p^\smallfrown s^\smallfrown b$
observes all the rules and the strategy $\gs$ and $t_p^\smallfrown s^\smallfrown b\subset t\}$. It is not difficult to
verify that the closed set $[S_p]$ is in the ideal $J_n$ since $i_p\mod n\notin c_n''[[S_p]]^2$. We also
have that $A\subset\bigcup_p[S_p]$: if some infinite binary sequence $x\in A$ fell out of all the sets $[S_p]$,
then a play of the game observing the strategy $\gs$ could be constructed such that the resulting sequence
is just $x$ and thus Adam won, contradicting the assumption that $\gs$ was a winning strategy for Eve. 

The previous two paragraphs together with the classical determinacy results of \cite{martin:borel} show that Borel sets
have the regularity property, and that ${\bf \gS}^1_n$ determinacy implies that ${\bf \gS}^1_n$ sets
have the regularity property. A standard trick described for example in \cite{kechris:bounded} can be used to reduce
the assumption to ${\bf\Delta}^1_n$ determinacy for the regularity of ${\bf \gS}^1_n$ sets, proving the Lemma.
\end{proof}

Thus the partial order $P_n$ is in forcing sense equivalent to the ordering of $n$-fat trees under inclusion.
It is now possible to give a detailed analysis of its forcing properties, using standard combinatorial methods.
However, the approach of the current paper is completely different. We shall need only the following consequence:

\begin{corollary}
\label{homogeneity}
The ideal $J_n$ is homogeneous. The forcing $P_n$ is homogeneous.
\end{corollary}

\noindent Here,

\begin{definition}
\cite{z:detandci} A $\gs$-ideal $I$ on $2^\gw$
is called \emph{homogeneous} if for every $I$-positive Borel set $B$ there is a function $f:2^\gw\to B$ such that
the preimages of $I$-small sets are $I$-small.
\end{definition}

Homogeneity of an ideal $I$ is a convenient way of securing the equalities $\cov(I)=\cov(I\restriction B)$
and $\non(I)=\non(I\restriction B)$ for an arbitrary positive Borel set $B$. 
If $X$ is a collection of $I$-small sets covering the set $B$ then the collection of their $I$-small preimages
covers the whole $2^\gw$; it has size $\leq |X|$. And if $Y\subset 2^\gw$ is an $I$-positive set, then its $f$-image is an $I$-positive
subset of the set $B$; it has size $\leq |Y|$. The key property of homogeneity is that it is frequently preserved
under the iterated Fubini powers of the ideal--\cite{z:detandci}.

\begin{proof}
Suppose that $B\subset 2^\gw$ is a $J_n$-positive Borel set, and let $i:2^\gw\to B$ be a continuous injection
respecting the partition $c_n$; such an injection exists by the previous lemma. Clearly, the preimages
of $J_n$-small sets must be $J_n$-small, and therefore the ideal $J_n$ is homogeneous. The function
$\pi$ defined by $\pi(C)=\pi''C$ for every $J_n$-positive Borel set $C\subset 2^\gw$, is an isomorphism
between the poset $P_n$ and $P_n$ below the condition $\pi(2^\gw)\subset B$. Note that Borel injective images of
Borel sets are Borel.
\end{proof}

\begin{corollary}
\label{definability}
The ideal $J_n$ is ${\bf\Pi^1_1}$ on ${\bf\gS^1_1}$.
\end{corollary}

\begin{proof}
Every $\gs$-ideal $K$ such that the poset $\cB(2^\gw)\setminus K$ ordered by inclusion is bounding
and contains a ${\bf\gS^1_1}$ collection of compact sets as a dense subset, is ${\bf\Pi^1_1}$ on
${\bf\gS^1_1}$--see the appendix of \cite{z:detandci}. A direct proof along the lines of
\cite{kechris:bounded} is also possible.
\end{proof}

Using Theorem 7.4 of \cite{z:detandci} and the previous lemmas it is also possible to conclude that
the countable support iteration of $P_n$ forcings is the optimal way to increase the invariant $\hm_n$.
While this is an interesting fact in itself, it plays no role in the proofs of the theorems from the introduction.

\begin{corollary}
(ZFC+LC) Suppose $n>0$ and $\fx$ is a tame invariant. If $\fx<\hm_n$ holds in some forcing extension, then
it holds in the countable support iterated $P_n$-extension.
\end{corollary}

As an aside, let us define the cardinal invariant $\hm_\gw$ in the following way. Let $\{a_n:n\in\gw\}$
be an arbitrary partition of $\gw$ into infinite sets, and for distinct sequences $x,y\in 2^\gw$
define $c_\gw(x,y)=n$ if the smallest number $m$ such that $x(m)\neq y(m)$ belongs to the set $a_n$.
Let $J_\gw$ be the $\gs$-ideal generated by the sets $X\subset 2^\gw$ such that $c_\gw''X\neq\gw$,
and define $\hm_\gw$ to be the covering number of this ideal. The invariant $\hm_\gw$ is independent 
of the initial choice of the partition $\{a_n:n\in\gw\}$: if $\{a'_n:n\in\gw\}$ is another such
a partition, $c'_\gw$ the associated function and $\hm'_\gw$ the associated invariant,
it is not difficult to find a continuous injection $i:2^\gw\to 2^\gw$ reducing the function
$c_\gw$ to $c'_\gw$ in the sense that $c'_\gw(i(x), i(y))=c_\gw(x,y)$ for every two distinct
sequences $x,y\in 2^\gw$. This shows that $\hm_\gw\leq\hm'_\gw$, and by symmetricity
$\hm_\gw=\hm'_\gw$. All the previous results all apply to the case
of the invariant $\hm_\gw$.   However, it is impossible to generalize the theorems
stated in the introduction to the case of $\hm_\gw$, since the results of the next section will not be applicable.

\section{The iterated Fubini powers of the ideals $J_n$}

The theorems stated in the introduction are proved using a careful analysis of the countable support
iteration of the forcings $P_n$. We will need to find an upper bound, in terms of the $\aleph$ function,
of the uniformities of the iterated Fubini powers of the ideals $J_n$. Recall:

\begin{definition}
\cite{z:detandci}
Suppose that $I$ is a $\gs$-ideal on the real line and $\ga\in\gw_1$ is an ordinal. The $\ga$-th iterated
Fubini power of the ideal $I$ is the $\gs$-ideal $I^\ga$ on $\R^\ga$ consisting of those sets $A\subset\R^\ga$
for which Adam has a winning strategy in the two person game $G(A)$ of length $\ga$. In $\gb$-th round
of the game $G(A)$ first Adam plays a set in the ideal $I$ and then Eve chooses a real not in the set.
Eve wins if the sequence of her answers belongs to the set $A$.
\end{definition}

But in fact,
in order to streamline the notation, we will not deal with the ideals $J_n$ directly. For every $n>0$ let
$K_n$ be the $\gs$-ideal on $(2^\gw)^n$ generated by the graphs of
Borel partial functions from $(2^\gw)^{n\setminus\{k\}}$ to $2^\gw$ for all $k\in n$. There
is a natural relationship between the ideals $K_n$ and $J_n$:
consider the bijection $g_n$ between $(2^\gw)^n$ and $2^\gw$ defined by $g_n(x_0, x_1, \dots
x_{n-1})=y$ where $y(m)=x_l(k)$ where $k$ is the integer part of $m/n$ and $l=m\mod n$. What is the
preimage of the ideal $J_n$? Let $X\subset 2^\gw$ be
a closed set such that $c_n''X$ misses at least one value, say $k\in n$. It is not difficult to see that
the preimage of the set $X$ under the bijection $g$ is a graph of a partial Borel function from
$(2^\gw)^{n\setminus\{k\}}$ to $2^\gw$: if $\vec x=\vec y$ are two $n$-tuples such that $g_n(\vec x), g_n(\vec y)
\in X$ and $\vec x(m)=\vec y(m)$ for all numbers $m\in n$ different from $k$, then necessarily $\vec x(k)=\vec y(k)$
as well, since otherwise $c_n(g_n(\vec x), g_n(\vec y))=k$, contradicting the choice of the set $X$.
Thus clearly $\{g^-1(X):X\in J_n\}\subset K_n$,
$\non(J_n)\leq\non(K_n)$ and for each countable ordinal $\ga$, $\non(J_n^\ga)\leq\non(K_n^\ga)$.
Thus it will be enough to find the uniformities of the ideals $K_n^\ga$.

\begin{lemma}
Let $n>0$ be a natural number and $\ga$ be a countable ordinal. Then $\non(K_n^\ga)\leq\aleph_{n+1}$.
\end{lemma}

\begin{proof}
The following is the key fact in the argument:

\begin{fact}
\label{clubguessing}
(Shelah) Let $\kappa>\gw_1$ be a regular cardinal. For every countable limit ordinal $\ga\in\gw_1$ there is a sequence
$\vec C=\langle C_\gd:\gd\in\kappa\rangle$ such that

\begin{enumerate}
\item $C_\gd\subset\gd$ is a closed set of ordertype $\leq\ga$
\item for every ordinal $\gd\in\kappa$ and every accumulation point $\gg\in C_\gd$, the set $C_\gg$ is just
$C_\gd\cap\gg$
\item for every closed unbounded set $E\subset\kappa$ there is an ordinal $\gd\in E$ such that the set $C_\gd$
is cofinal in $\gd$ of ordertype $\ga$ and it is a subset of $E$.
\end{enumerate}
\end{fact}

This Fact was announced in \cite{shelah:nonstructure}, page 136, remark 2.14A. The proof
is in \cite{math.LO/9906022}, available from the Mathematics ArXiv. 
Now fix a number $n>0$ and a countable ordinal $\ga\in\gw_1$, without loss
assuming that $\ga$ is limit. For every number $m\in n$ use the Fact to choose a club guessing sequence
$\vec C_m=\langle C^m_\gd:\gd\in\gw_{m+2}\rangle$ on $\gw_{m+2}$. We may certainly assume that $\fc\geq\aleph_{n+1}$
and so we can choose a sequence $\vec s=\langle s_\gd:\gd\in\gw_{n+1}\rangle$ of pairwise distinct reals--elements
of $2^\gw$. For every $n$-tuple $\vec\gd\in\prod_{m\in n}\gw_{m+2}$ such that the sets $C^m_{\vec\gd(m)}$ have
ordertype $\ga$, let $\langle\vec\gd(m)(\gb):\gb\in\ga\rangle$ enumerate the nonaccumulation points of these sets
in the increasing order, and let $\vec r_{\vec \gd}\in(((2^\gw)^n)^\ga)$ be the $\alpha$-sequence whose
$\gb$-th element is the point $\langle s_{\vec\gd(m)(\gb)}:m\in n\rangle$ in the space $(2^\gw)^n$. It will be enough to
show that the set $\{\vec r_{\vec\gd}:\vec\gd\in\prod_{m\in n}\gw_{m+2}\}$ is $J^\ga_n$-positive.

To prove this, for every Adam's strategy $\gs$ in the transfinite game we need to find an $n$-tuple $\vec\gd$ such that
the sequence $\vec r_{\vec \gd}$ is a legal counterplay against the strategy $\gs$. So fix the strategy $\gs$
and by downward induction on $m\in n$ find

\begin{enumerate}
\item a continuous increasing $\in$-tower $\vec M_m=\langle M^m_\gd:\gd\in\gw_{m+2}\rangle$ of elementary submodels
of a large enough structure, each of them of size $\aleph_{m+1}$ and such that $M^m_\gd\cap\gw_{m+2}\in\gw_{m+2}$.
In particular, we require $\vec C_k:k\in n+1, \vec s, \gs\in M^m_0$. 
Let $E_m=\{\gd\in\gw_{m+2}:M^m_\gd\cap\gw_{m+2}=\gd\}$; this
is a closed unbounded subset of $\gw_{m+2}$.
\item an ordinal $\gd_m\in\gw_{m+2}$ such that the set $C^m_{\gd_m}\subset\gd_m$ is cofinal of ordertype $\ga$
and it is a subset of the club $E_m$. For every two numbers $k\in m\in n$ we demand that $\gd_m\in M^k_0$.
\end{enumerate}

Let $\vec\gd=\langle\gd_m:m\in n\rangle$. We claim that the sequence $\vec r_{\vec\gd}$ is the desired legal
counterplay against the strategy $\gs$. So look at an arbitrary round $\gb\in\ga$ and suppose that the
sequence does constitute a legal counterplay up to this point. What happens at round $\gb$? 

The important observation is that for every integer $k\in n$,
the play up to this round is in the model  $M^k_{\vec\gd(k)(\gb)}$ since it is defined
from objects that belong to this model. In particular, one of the parameters in the definition is the set
$C^k_{\vec\gd(k)}\cap\vec\gd(k)(\gb)$ which is in the model by the coherence requirement (2) in Fact~\ref{clubguessing}.

The strategy $\gs$
now commands Adam to play partial Borel functions $\{f_{mk}:k\in n,m\in\gw\}$ such that $f_{mk}$ is a function
from $(2^\gw)^{n\setminus\{k\}}$ to $2^\gw$. We must show that the point $\vec r_{\vec \gd}(\gb)$
is not contained in the graph of any of these functions. So choose integers $m\in\gw$ and $k\in n$.
Consider the set $Y=\{f_{mk}(\vec u,\vec v):\vec u\in(2^\gw)^k$ is a sequence all of whose entries 
are on the $\vec s$-sequence,
indexed by ordinals smaller than $\gw_{k+1}$ and $\vec v\in(2^\gw)^{n\setminus(k+1)}$ is a sequence
all of whose entries are on the $\vec s$-sequence, indexed by ordinals 
in the set $\bigcup_{k+2\in l\in n+2}C^l_{\vec\gd(l)}\}$.
This set is of size $<\gw_{k+2}$ and it belongs to the model $M^k_{\vec\gd(k)(\gb)}$.
Thus $Y\subset M^k_{\vec\gd(k)(\gb)}$, in particular $s_{\vec\gd(k)(\gb)}\notin Y$, which by the definition of
the set $Y$ means that the point $\vec r_{\vec\gd}(\gb)=\langle s_{\vec\gd(l)(\gb)}:l\in n\rangle$ is not on the
graph of the function $f_{mk}$ as desired. The lemma follows.
\end{proof}

\begin{corollary}
\label{bab}
Let $n>0$ be a natural number and let $\ga$ be a countable ordinal. $\non(J_n^\ga)\leq\aleph_{n+1}$.
\end{corollary}

In order to analyse the countable support iteration in which the forcings $P_n$ alternate, we need to change the approach
a little bit:

\begin{definition}
Suppose that $\vec I$ is an $\gw$-sequence of $\gs$-ideals on the reals and $\ga$ is a countable ordinal. The
ideal $\vec I^\ga$ on $\R^{\gw\cdot\ga}$ consists of those sets $A\subset\R^{\gw\cdot\ga}$ for which Adam has a
winning strategy in the transfinite two person game $G(A)$ of length $\gw\cdot\ga$. In the $\gw\cdot\gb+n$-th
round of this game, Adam plays a set in the ideal $\vec I(n)$ and Eve responds with a real which is not in this set.
Eve wins if the sequence of her answers belongs to the set $A$.
\end{definition}

Let $\vec K$ be the $\gw$-sequence of ideals given by $\vec K(n)=K_{n+1}$.

\begin{lemma}
Let $\ga$ be a countable ordinal. Then $\non(\vec K^\ga)\leq\aleph_{\gw_2+1}$.
\end{lemma}

\begin{proof}
Let $\ga$ be an arbitrary countable ordinal. First we need to fix several objects whose existence is provable in ZFC.

\begin{enumerate}
\item An increasing sequence $\vec \kappa=\langle\kappa_\gd:\gd\in\gw_2\rangle$ of regular cardinals below $\aleph_{\gw_2}$
such that the true cofinality of their product modulo the bounded ideal on $\gw_2$ is $\oox$,
from \cite{shelah:pcf}, Chapter II, Theorem 1.5. This means that there is a sequence 
$\vec h=\langle h_\gg:\gg\in\oox\rangle$ of functions in $\prod_\gd\kappa_\gd$ which is increasing and cofinal
in the modulo bounded ordering. Fix such a sequence.
\item A club guessing sequence $\vec C=\langle C_\gg:\gg\in\oox\rangle$ from Fact~\ref{clubguessing}. The sequence will guess closed
unbounded subsets of $\oox$ by segments of length $\gw\cdot\gw\cdot\ga$. For every ordinal $\gg\in\oox$
let $C(\gg, \gb)$ denote the $\gb$-th nonaccumulation point of the set $C_\gg$.
\item A club guessing sequence $\vec D=\langle D_\gd:\gd\in\gw_2\rangle$ from Fact~\ref{clubguessing}. The sequence will guess closed unbounded
subsets of $\gw_2$ by segments of length $\gw\cdot\gw\cdot\ga$ again, with similar notational convention
for as in the previous
item, using the symbol $D(\gd,\gb)$ for the $\gb$-th nonaccumulation point of the set $D_\gd$.
\item Without harm we may assume that $\fc>\aleph_{\gw_2}$. So let us fix a sequence $\vec s=\langle s_\gg:\gg\in\oox\rangle$
of pairwise distinct reals (elements of $2^\gw$).
\end{enumerate}

Now suppose that $\gg\in\oox$ and $\gd\in\gw_2$ are ordinals such that the ordertypes of the sets $C_\gg$ and $D_\gd$
are both $\gw\cdot\gw\cdot\ga$. Define a $\gw\cdot\ga$ sequence $\vec r_{\gg\gd}$ by setting its $\gw\cdot\gb+n$-th
element to be the point in the space $(2^\gw)^{n+1}$ whose $k$-th coordinate for every number $k\in n+1$ is
the real number on the $s$ sequence indexed by the ordinal $h_{C(\gg, \gw\cdot(\gw\cdot\gb+n)+k)}(D(\gd,\gw\cdot
(\gw\cdot\gb+n)+n-k))$. We will show that the collection $\{\vec r_{\gg\gd}:\gg\in\oox,\gd\in \gw_2\}$ is
$\vec K^\ga$-positive, proving the lemma. This means that for every Adam's strategy $\gs$ in the transfinite game
we must find ordinals $\gg\in\oox$ and $\gd\in\gw_2$ such that the sequence $\vec r_{\gg\gd}$ is a legal counterplay
against the strategy.

Fix a continuous increasing $\in$-tower $\langle M_\gg:\gg\in\oox\rangle$ of elementary submodels of large enough structure,
each of them of size $\aleph_{\gw_2}$ and such that $M_\gg\cap\oox\in\oox$. In particular, 
$\vec \kappa, \vec h, \vec C, \vec D, \gs\in M_0$. Let $E=\{\gg\in\oox:M_\gg\cap\oox\in\oox\}$.
Since this is a closed unbounded subset of $\oox$, there must be an ordinal $\gg$ such that the set $C_\gg\subset\gg$
is cofinal of ordertype $\gw\cdot\gw\cdot\ga$ and $C_\gg\subset E$.

Also, fix a continuous increasing $\in$-tower $\langle N_\gd:\gd\in\gw_2\rangle$ of elementary submodels of
large enough structure, each of them of size $\aleph_1$ and such that $N_\gd\cap\gw_2\in\gw_2$. In particular,
$\vec \kappa,\vec h, \vec C, \vec D, E,\gs,\gg\in N_0$. 
Let $F=\{\gd\in\gw_2:N_\gd\cap\gw_2=\gd\}$. Since this is a closed unbounded set, there
must be an ordinal $\gd\in\gw_2$ such that the set $D_\gd\subset\gd$ is cofinal of ordertype $\gw\cdot\gw\cdot\ga$
and $D_\gd\subset F$.

We claim that the sequence $\vec r_{\gg\gd}$ is a legal counterplay against the strategy $\gs$. Well, consider
the situation at round $\gw\cdot\gb+n$ for some ordinal $\gb\in\ga$ and number $n\in\gw$. Suppose that
up to this point, the sequence constituted a legal partial play; we want to see that it will provide
a legal answer even in this round. The strategy $\gs$ commands Adam to play partial Borel functions
$\{f_{mk}:k\in n+1,m\in\gw\}$ such that for each $k\in n+1$ and each $m\in\gw$ the function
maps $(2^\gw)^{n+1\setminus\{k\}}$ to $(2^{\gw})$. We must show that the $n+1$-tuple $\vec t=\vec r_{\gg\gd}
(\gw\cdot\gb+n)\in(2^\gw)^{n+1}$ is not on the graph of any of these functions, that is $\vec t(k)\neq
f_{mk}(\vec t\restriction(n+1\setminus\{k\}))$. 

To this end, fix integers $k\in n+1$ and $m\in\gw$. Define a function $g\in\prod_\eta\kappa_\eta$ by
letting $g(\eta)$ to be the supremum of the set $\{\xi\in\kappa_\eta:s_\xi=f_{mk}(\vec u,\vec v),$ where
$\vec u\in(2^\gw)^k$ is a sequence all of whose entries are on the $\vec s$-sequence, and are indexed by
ordinals $<\sup\{\kappa_{\eta'}:\eta'\in\eta\}$, and $\vec v\in(2^\gw)^{(n+1)\setminus(k+1)}$ is a sequence
all of whose entries are on the $\vec s$-sequence and are indexed in the ordinals in the range of the functions
$\{h_{\gg'}:\gg'\in C_\gg\cap C(\gg ,\gw\cdot(\gw\cdot\gb+n)+k)\}\}$. 
It is immediate that this set has size $<\kappa_\eta$ and
so the function $g$ is well-defined. There are two important points.

\begin{enumerate}
\item $g\in M_{C(\gg ,\gw\cdot(\gw\cdot\gb+n)+k)}$. This so happens because the function is defined from objects
contained in the model, in particular from the set $C_\gg\cap\gg(\gw\cdot(\gw\cdot\gb+n)+k)$ which belongs to
the model by the coherence of the $C$-sequence.
\item $g\in N_{D(\gd,\gw\cdot(\gw\cdot\gb+n)+k)}$ by the same reason as in the previous item, this time using
the coherence of the $D$-sequence.
\end{enumerate}

By the first point, the function $g\in\prod_\gd\kappa_\gd$ is dominated by the function $h_{C(\gg,\gw\cdot(\gw\cdot\gb+n)+k)}$ from
some ordinal on. By the second point, this ordinal must be smaller than $D(\gd,\gw\cdot(\gw\cdot\gb+n)+k)$. By the definition
of the function $g$ and the sequence $\vec t=\vec r_{\gg\gd}
(\gw\cdot\gb+n)$ then, it must be the case that $\vec t(k)\neq
f_{mk}(\vec t\restriction(n+1\setminus\{k\}))$ as desired. The lemma follows.
\end{proof}

So we have the following. Let $\vec J$ be the $\gw$-sequence of ideals given by $\vec J(n)=J_{n+1}$.

\begin{corollary}
\label{dad}
Let $\ga$ be a countable ordinal. $\non(\vec J^\ga)\leq\aleph_{\gw_2}+1$.
\end{corollary} 

Let $K_\gw$ be the $\gs$-ideal on $(2^\gw)^\gw$ $\gs$-generated by the graphs of partial Borel functions from
$(2^\gw)^{\gw\setminus\{n\}}$ to $2^\gw$, for all numbers $n\in\gw$. The following lemma is an observation
complementary to the previous results in this section. It shows that the ideal $K_\gw$ is different from the
ideals $K_n$ in that its uniformity can be arbitrarily large:

\begin{lemma}
$\mathfrak{p}\leq\non(K_\gw)$.
\end{lemma}

\begin{proof}
Suppose that $A\subset(2^\gw)^\gw$ is a set of size $<\mathfrak{p}$. We must produce partial Borel functions
$\{f_{n}:n\in\gw\}$ such that $f_n:(2^\gw)^{\gw\setminus\{n\}}\to 2^\gw$ and every point in the set $A$
is on the graph of one of them. First, use the Axiom of Choice to find a set $B\subset A$ so that

\begin{enumerate}
\item for all sequences $\vec x,\vec y\in B$ if $\vec x(m)=\vec y(m)$ for all but finitely many integers $m$ then $\vec x=
\vec y$
\item $B\subset A$ is a maximal set with the previous property.
\end{enumerate}

For each number $n$ let $g_n$ be a partial function from $(2^\gw)^{\gw\setminus\{n\}}$ to $2^\gw$ defined by
$g_n(\vec x)=z$ if the sequence $\vec y_0=\vec x\cup\{\langle n, z\rangle\}$ is in the set $A$ and for some sequence
$\vec y_1$ in the set $B$, $\vec y_0(m)=\vec y_1(m)$ for all numbers $m\geq n$. By the property (1) of the set $B$,
this formula does correctly define the functions $g_n$, and by the property (2), every point in the set $A$ is
on the graph of all but finitely many of these functions. Now there is a general fact, proved for example in \cite{z:detandci}
Appendix B, that every partial function of size $<\mathfrak{p}$ between Polish spaces is a subset
of a Borel function. Thus for every number $n$ there is a Borel function $f_n$ such that $g_n\subset f_n$;
clearly every point in the set $A$ is on the graph of all but finitely many of the functions $f_n$ as desired.

To prove the abovementioned general fact, consider for simplicity a partial function $g:2^\gw\to 2^\gw$.
By a standard almost disjoint coding argument find a $\gs$-centered forcing adding subsets $\dot a_i:i\in\gw$
of $2^{<\gw}$ such that $x\cap\dot a_i$ is finite if and only $g(x)(i)=0$ for every infinite binary sequence
$x$ in the ground model, viewed as the set of all its initial segments and so a subset of $2^{<\gw}$.
By the small size of the function $g$ and Bell's theorem there will be sets $a_i:i\in\gw$ already
in the ground model such that the Borel function $f:2^\gw\to 2^\gw$ defined by $f(x)(i)=0$ iff $x\cap a_i$
is finite, extends the function $g$ as desired.
\end{proof}

\section{The duality theorems}

The proof of Theorem~\ref{t2} follows the argument for the duality theorem in the Applications section
of \cite{z:detandci}. Suppose that $n>0$ is a natural number and $I$ is an analytic 
ideal on the reals such that ZFC proves $\cov(I)\geq\hm_n$. We want to argue that 
$\non(I)\leq\aleph_{n+1}$. For this, it is necessary to analyse the countable support iteration of the $P_n$
forcing. The following Facts use the terminology and arguments from \cite{z:detandci}. Recall:

\begin{definition}
\cite{z:detandci}
Suppose $K$ is a $\gs$-ideal on the reals and $\ga\in\gw_1$ is an ordinal. A set $B\subset\R^\ga$
is $K$-perfect if the tree $T\subset\R^{\leq\ga}$ of all initial segments of the sequences in the set $B$
has the following two properties. It is $K$-positively
branching, meaning that for every node $t\in T$ the set $\{r\in\R:t^\smallfrown r\in T\}$ is not in the ideal $K$.
And it is $\gs$-closed, meaning that for every sequence $t_0\subset t_1\subset\dots$ of sequences in the tree
$T$ we have $\bigcup_{n\in\gw}t_n\in T$.
\end{definition}

\begin{fact}
\label{eee}
For every countable ordinal $\ga\in\gw_1$ and every analytic set $A\subset(2^\gw)^\ga$, either
the set $A$ contains a Borel $J_n$-perfect subset, or it belongs to the ideal $J_n^\ga$. Under
the assumption of the existence of $\gw$ many Woodin cardinals, this extends to all projective sets.
\end{fact}

This is an immediate consequence of Lemmas~\ref{l1} and~\ref{l5}, Corollary~\ref{definability} 
and the work of \cite{z:detandci}.

\begin{fact}
For every countable ordinal $\ga\in\gw_1$ the ideal $J_n^\ga$ is homogeneous.
\end{fact}

This happens because the ideal $J_n$ is homogeneous. It means in particular that for every Borel $J_n$-perfect
set $B\subset(2^\gw)^\ga$ the invariant $\non(J_n^\ga)$ is equal to $\non(J_n^\ga\restriction B)$, the smallest size of 
a $J_n^\ga$-positive subset of $B$.

\begin{fact}
There is a countable ordinal $\ga\in\gw_1$, a Borel $J_n$-perfect subset $B\subset(2^\gw)^\ga$ and a Borel
function $f:B\to\R$ such that $f$-preimages of $I$-small sets are $J_n^\ga$-small.
\end{fact}

This is true because after the countable support iteration of the $P_n$-forcing of length $\fc^+$, it is still the
case that $\cov(I)\geq\hm_n$ and so there must be a name for a real which falls out of all ground model coded
$I$-small sets. This name must be represented by a Borel set and a function as in the above Fact--\cite{z:detandci}.

Now fix $\ga, B$ and $f$ from the previous Fact. By Fact~\ref{homogeneity} and Corollary~\ref{bab}, there is
a $J^\ga_n$-positive subset $A\subset B$ of size $\leq\aleph_{n+1}$. By the properties of the function
$f$, the set $f''A$ must be an $I$-positive set of size $\leq\aleph_{n+1}$, and Theorem~\ref{t2} follows.

The argument for Theorem~\ref{t3} is similar,
using the analysis of the countable support iteration in which the forcings $P_n$ alternate and Corollary~\ref{dad}.
Theorem~\ref{t1} is just a special case of Theorem~\ref{t2}.

\def\germ{\frak} \def\scr{\cal} \ifx\documentclass\undefinedcs
  \def\bf{\fam\bffam\tenbf}\def\rm{\fam0\tenrm}\fi 
  \def\defaultdefine#1#2{\expandafter\ifx\csname#1\endcsname\relax
  \expandafter\def\csname#1\endcsname{#2}\fi} \defaultdefine{Bbb}{\bf}
  \defaultdefine{frak}{\bf} \defaultdefine{mathfrak}{\frak}
  \defaultdefine{mathbb}{\bf} \defaultdefine{mathcal}{\cal}
  \defaultdefine{beth}{BETH}\defaultdefine{cal}{\bf} \def\bbfI{{\Bbb I}}
  \def\mbox{\hbox} \def\text{\hbox} \def\om{\omega} \def\Cal#1{{\bf #1}}
  \def\pcf{pcf} \defaultdefine{cf}{cf} \defaultdefine{reals}{{\Bbb R}}
  \defaultdefine{real}{{\Bbb R}} \def\restriction{{|}} \def\club{CLUB}
  \def\w{\omega} \def\exist{\exists} \def\se{{\germ se}} \def\bb{{\bf b}}
  \def\equivalence{\equiv} \let\lt< \let\gt>

\end{document}